



\hsize=12cm
\vsize=19cm
\hoffset=2cm
\voffset=2cm
\footline={\hss{\vbox to 2cm{\vfil\hbox{\rm\folio}}}\hss}

\input amssym.def
\input amssym.tex

\font\csc=cmcsc10
\font\title=cmr12 at 14pt
\font\stitle=cmsl12 at 14pt

\font\teneusm=eusm10    
\font\seveneusm=eusm7  
\font\fiveeusm=eusm5    
\newfam\eusmfam
\def\eusm{\fam\eusmfam\teneusm}
\textfont\eusmfam=\teneusm 
\scriptfont\eusmfam=\seveneusm
\scriptscriptfont\eusmfam=\fiveeusm

\def\varGamma{{\mit \Gamma}}
\def\Re{{\rm Re}\,}
\def\Im{{\rm Im}\,}
\def\txt#1{{\textstyle{#1}}}
\def\scr#1{{\scriptstyle{#1}}}
\def\r#1{{\rm #1}}
\def\B#1{{\Bbb #1}}
\def\e#1{{\eusm #1}}

\def\scr#1{{\scriptstyle #1}}

\centerline{\title A Functional Equation}
\medskip
\centerline{\title 
for the Spectral Fourth Moment}
\medskip
\centerline{\title of Modular Hecke {\stitle L}\thinspace-Functions}
\vskip 1cm
\centerline{\csc By Y. Motohashi}
\vskip 1cm
\par
\noindent
Albeit essential corrections are
required both in his claim and in his argument, N.V. Kuznetsov
observed in [2] a highly interesting transformation formula for spectral
sums of products of four values of modular Hecke $L$-functions. A
complete proof of a corrected version of the formula  was later
supplied by the present author [3], which has, however, remained
unpublished for more than a decade, except for a limited distribution. 
The aim of the present paper is to reproduce [3] with some
sophistications, retaining the original style as much
as possible.
\smallskip
\noindent
{\csc Convention}: Notations become available at their first
appearances and will continue to be effective thereafter.
\bigskip
\noindent 
{\bf 1.}\quad We shall first introduce basic notions; for
details we refer to [4]. Thus, let
$\varGamma$ be the full modular group $\r{PSL}_2(\B{Z})$, and
$L^2(\varGamma\backslash\B{H})$ the Hilbert space composed of all
$\varGamma$-automorphic functions on 
$\B{H}=\{x+iy: y>0\}$, which are square integrable over
$\varGamma\backslash\B{H}$ against
the hyperbolic measure. A real analytic modular cusp form or a Maass
form is an element in
$L^2(\varGamma\backslash\B{H})$, which is an eigenfunction of the
hyperbolic Laplacian $\Delta=-y^2((\partial/\partial
x)^2+(\partial/\partial y)^2)$. The subspace spanned by all Maass forms
has a complete orthonormal system
$\{\psi_j: j=1,2,\cdots\}$ such that
$\Delta\psi_j=({1\over4}+\kappa_j^2)\psi_j$ with 
$0<\kappa_1\le\kappa_2\le\cdots$, and
$$
\psi_j(x+iy)=\sqrt{y}\sum_{n\ne0}\varrho_j(n)
K_{i\kappa_j}(2\pi|n|y)e(nx),\quad x+iy\in\B{H},\eqno(1.1)
$$
where $K_\nu$ is the $K$-Bessel function of order $\nu$, and
$e(x)=\exp(2\pi ix)$. The $\varrho_j(n)$ are called the
Fourier coefficients of $\psi_j$. We may assume as usual
that  $\psi_j$  are simultaneous
eigenfunctions of all Hecke operators with corresponding
eigenvalues
$t_j(n)\in\B{R}$; that is, for each integer $n>0$
$$
{1\over\sqrt{n}}\sum_{ad=n}\sum_{b\bmod d}
\psi_j((az+b)/d)
=t_j(n)\psi_j(z),\quad z\in\B{H},\eqno(1.2)
$$
and that
$$
\psi_j(-\overline{z})=\epsilon_j\psi_j(z),
\quad\epsilon_j=\pm1.\eqno(1.3)
$$
We have, uniformly in $j$,
$$t_j(n)\ll
n^{{1\over4}+\varepsilon}\eqno(1.4)
$$
for any fixed $\varepsilon>0$.
The Hecke $L$-function associated with
$\psi_j$ is defined by
$$
H_j(s)=\sum_{n=1}^\infty t_j(n)n^{-s},\quad \Re s>1.\eqno(1.5)
$$
This continues to an entire function satisfying a Riemannian functional
equation; in consequence, $H_j(s)$
is of polynomial growth with respect to both $\kappa_j$ and $s$ in any
fixed vertical strip. 
\par
We need also to define holomorphic modular cusp forms, and
corresponding Hecke $L$-functions. Thus, if $\psi$ is
holomorphic over $\B{H}$, and $\psi(z)(dz)^k$ with an
integer $k>0$ is $\varGamma$ invariant, then $\psi$ is 
a holomorphic modular form of weight $2k$. If $\psi$ vanishes
at the point infinity, then it is called a cusp form. The set
composed of all such functions is a finite dimensional Hilbert
space. Let $\{\psi_{j,k}:1\le j\le\vartheta(k)\}$
be an orthonormal basis. The Fourier coefficients $\varrho_{j,k}(n)$ of
$\psi_{j,k}$ is defined by the expansion
$$
\psi_{j,k}(z)=\sum_{n=1}^\infty
n^{k-{1\over2}}\varrho_{j,k}(n)e(nz),\quad z\in\B{H}.\eqno(1.6)
$$
We may assume again that $\psi_{j,k}$ are simultaneous
eigenfunctions of all Hecke operators, so that there exist
real numbers $t_{j,k}(n)$ such that
$$
{1\over\sqrt{n}}\sum_{ad=n}(a/d)^k
\sum_{b\bmod d}\psi_{j,k}((az+b)/d)
=t_{j,k}(n)\psi_{j,k}(z).\eqno(1.7)
$$
Corresponding to $(1.4)$,
we have, uniformly in $j,k$,
$$
t_{j,k}(n)\ll n^{{1\over4}+\varepsilon}.\eqno(1.8)
$$
Also the Hecke $L$-function attached to $\psi_{j,k}$ is defined
by
$$
H_{j,k}(s)=\sum_{n=1}^\infty t_{j,k}(n)n^{-s},\quad
\Re s>1,\eqno(1.9)
$$
which continues to an entire function of polynomial growth
with respect to both $k$ and $s$ in any fixed vertical strip.
\bigskip
\noindent{\bf 2.}\quad Now, let $h(r)$ be even,
entire, and
$$
h(r)\ll\exp(-c|r|^2).\eqno(2.1)
$$
This assumption is made for the sake of simplicity; in fact, we can
relax it considerably. For instance, our argument works well, provided 
that $h(r)$ is regular and $\ll (|r|+1)^{-A}$ for $|\Im r|<A$
with a sufficiently large
$A$. 
\smallskip
In the region of absolute convergence, i.e., 
$$
z=(z_1,z_2,z_3,z_4)\in\B{C}^4,\quad\Re z_j>1\quad (1\le j\le4),
\eqno(2.2)
$$
let 
$$
\e{H}_\pm(z;h)=\sum_{j=1}^\infty\alpha_j\epsilon_j^{{1\over2}(1\mp1)}
H_j(z_1)H_j(z_2)
H_j(z_3)H_j(z_4)h(\kappa_j)\eqno(2.3)
$$
with $\alpha_j=|\varrho_j(1)|^2/\cosh(\pi\kappa_j)$, and
$$
\eqalignno{
&\e{C}(z;h)={1\over\pi}\int_{-\infty}^\infty\zeta(z_1+ir)\zeta(z_2+ir)
\zeta(z_3+ir)\zeta(z_4+ir)\cr
&\times\zeta(z_1-ir)\zeta(z_2-ir)
\zeta(z_3-ir)\zeta(z_4-ir)
{h(r)\over \zeta(1+2ir)\zeta(1-2ir)}dr&(2.4)
}
$$
with $\zeta$ the Riemann zeta-function. 
Put
$$
\e{G}(z;h)=\left(\e{H}_++\e{C}\right)(z;h).
\eqno(2.5)
$$
Obviously $\e{H}_+(z;h)$ is
entire over $\B{C}^4$. On the other hand, $\e{C}(z;h)$ is meromorphic,
which can be seen by shifting the contour vertically. Thus,
$\e{G}(z;h)$ is meromorphic over $\B{C}^4$. Interesting situations
occur when $z$ is in the set
$$
\e{X}=\left(\txt{1\over2}+i\B{R}\right)^4.\eqno(2.6)
$$
or in its vicinity. Then we have
$$
\e{G}(z;h)=\left(\e{H}_++\e{C}_0+\e{R}_0\right)(z;h),\eqno(2.7)
$$
where $\e{C}_0$ has the same integral representation as $\e{C}$ but
with the present specification for $z$, and
$$
\eqalignno{
&\e{R}_0(z;h)\cr
&=4\zeta(z_2-z_1+1)\zeta(z_3-z_1+1)\zeta(z_4-z_1+1)\cr
&\phantom{4}\times\zeta(z_2+z_1-1)\zeta(z_3+z_1-1)\zeta(z_4+z_1-1)
h(i(z_1-1))/\zeta(3-2z_1)\cr
&+4\zeta(z_1-z_2+1)\zeta(z_3-z_2+1)\zeta(z_4-z_2+1)\cr
&\phantom{4}\times\zeta(z_1+z_2-1)\zeta(z_3+z_2-1)\zeta(z_4+z_2-1)
h(i(z_2-1))/\zeta(3-2z_2)\cr
&+4\zeta(z_1-z_3+1)\zeta(z_2-z_3+1)\zeta(z_4-z_3+1)\cr
&\phantom{4}\times\zeta(z_1+z_3-1)\zeta(z_2+z_3-1)\zeta(z_4+z_3-1)
h(i(z_3-1))/\zeta(3-2z_3)\cr
&+4\zeta(z_1-z_4+1)\zeta(z_2-z_4+1)\zeta(z_3-z_4+1)\cr
&\phantom{4}\times\zeta(z_1+z_4-1)\zeta(z_2+z_4-1)\zeta(z_3+z_4-1)
h(i(z_4-1))/\zeta(3-2z_4)\qquad&(2.8)
}
$$
(cf. pp.\ 117--118 of [4]).
It can be shown, with an elementary argument, that $\e{R}_0$ is finite
on $\e{X}$.
\smallskip
\noindent
{\csc Remark}: Both [2] and [3] treat $\e{G}$. One may include
$\e{H}_-$ into discussion to make the subject more symmetric: The sum
$$
\e{H}_+(z;h_+)+\e{H}_-(z;h_-)+\e{C}(z;h_++h_-)\eqno(2.9)
$$
could be chosen to start with. To attain a fuller symmetry, the sum
$\e{H}_0$, defined by $(5.20)$ below, is also to be included. 
We note that dealing with a spectral sum
involving $H_j({1\over2})$  one may exploit the fact that 
$H_j({1\over2})=\epsilon_jH({1\over2})$ (see
Section 3.3 of [4]), and that $\e{H}_-$ is much
less problematic than $\e{H}_+$.
Hence, sometimes, $\e{H}_-$ can be a better choice than
$\e{H}_+$ to work with.
These refinements are nevertheless omitted in our present
discussion, because of their peripheral nature.
\smallskip
Returning to the domain $(2.2)$, we have
$$
\eqalignno{
&\e{H}_+(z;h)=\zeta(z_1+z_2)\zeta(z_3+z_4)\cr
&\times\sum_{m,n\ge1}{\sigma_{z_1-z_2}(m)\sigma_{z_3-z_4}(n)\over
m^{z_1}n^{z_3}}\sum_{j\ge1}\alpha_j
t_j(m)t_j(n)h(\kappa_j),&(2.10)
}
$$
by the multiplicativity of Hecke eigenvalues (see $(3.2.7)$ of [4]).
To the inner sum we apply the Spectral-Kloosterman sum formula due to
Bruggeman and Kuznetsov (Theorems 2.2 of [4]). We have
$$
\eqalignno{
&\e{G}(z;h)={1\over\pi^2\zeta(z_1+z_2+z_3+z_4)}
\cr
&\times\zeta(z_2+z_3)\zeta(z_2+z_4)\zeta(z_3+z_4)
\int_{-\infty}^\infty rh(r)\tanh\pi r\,dr\cr
&+\zeta(z_1+z_2)\zeta(z_3+z_4)
\sum_{m,n\ge1}{\sigma_{z_1-z_2}(m)\sigma_{z_3-z_4}(n)\over
m^{z_1}n^{z_3}}\e{K}(m,n;\varphi).&(2.11)
}
$$
Here 
$$
\e{K}(a,b;\varphi)=\sum_{c\ge1}{1\over c}S(a,b;c)
\varphi\left({4\pi\over l}\sqrt{|ab|}\right),\eqno(2.12)
$$
where $S(m,n;l)$ is the Kloosterman sum, and, for $x>0$,
$$
\varphi(x)={1\over\pi i}\int_{-\infty}^\infty 
(J_{-2ir}(x)-J_{2ir}(x)){rh(r)\over
\cosh \pi r}dr,\eqno(2.13)
$$
with $J_\nu$ the $J$-Bessel function of order $\nu$.
Considering the Mellin transform of $\varphi$, we have
$$
\varphi(x)={1\over\pi^2}\int_{(\alpha)}h^\ast
(s)(x/2)^{-2s}ds,\quad -{1\over2}<\alpha<0,
\eqno(2.14)
$$
where $(\alpha)$ is the vertical line $\Re s=\alpha$, and
$$
h^\ast(s)=\int_{-\infty}^\infty{\Gamma(s+ir)\over\Gamma(1-s+ir)}
{rh(r)\over\cosh \pi r}dr.\eqno(2.15)
$$
A downward shift of the contour shows that $h^*(s)$ is
regular for $\Re s>-{1\over2}$, and there
$$
h^*(s)\ll (|s|+1)^{2\Re s-1},\eqno(2.16)
$$
which is in fact the best possible.
However, this bound is not adequate to carry out subsequent
transformations of $\e{G}$, a grave fact which is overlooked in
[2]. 
\par
On the other hand the corresponding transform of $h$ attached to
$\e{H}_-$ is
$$
\int_{-\infty}^\infty{\Gamma(s+ir)\over\Gamma(1-s+ir)}
{rh(r)\over\cos\pi s}dr.\eqno(2.17)
$$
This decays exponentially, which makes $\e{H}_-$ much
easier to treat than $\e{H}_+$. See also p.\ 113 of [4].
\bigskip
\noindent
{\bf 3.}\quad To overcome this difficulty with $h^*$, we appeal to a
device that has stemmed from an idea of R. Murty, which exploits the
fact that the full modular group lacks holomorphic cusp forms of low
weights (see pp.\ 442--443 of [5] and Concluding Remark below).
\smallskip
Put
$$
h^{**}(s)=h^*(s)+{1\over2}\pi i\sum_{\scr{k=1}\atop\scr{k\ne6}}^7
\gamma_k{\Gamma(k-{1\over2}+s)\over\Gamma(k+{1\over2}-s)},\eqno(3.1)
$$
where for $1\le k\le3$
$$
\gamma_k=(-1)^k{2\over\pi}(2k-1)h(i(k-\txt{1\over2})),
\eqno(3.2)
$$
and other $\gamma_k$ are to satisfy, for $\nu=0,1,2$,
$$
\sum_{\scr{k=1}\atop\scr{k\ne6}}^7(-1)^k
\gamma_k(k-\txt{1\over2})^{2\nu}
=(-1)^\nu{2\over\pi}\int_{-\infty}^\infty r^{2\nu+1}
h(r)\tanh \pi r\,dr.\eqno(3.3)
$$
Obviously the $\gamma_k$ are fixed uniquely. 
With this, $h^{**}(s)$ is
regular over $\B{C}$ except for the simple
poles at $s=-(k-{1\over2})$ $(k=4,5,\ldots)$. Moreover, as $s$
tends to infinity in any fixed vertical strip, we have
$$
h^{**}(s)\ll |s|^{2\Re s-4}.\eqno(3.4)
$$
\par
The first assertion is obtained by shifting the contour in $(2.15)$
to $\Im r=-3$. The second will  be proved later for the
sake of completeness. We remark that the 
modification $(3.1)$ does not
cause much changes in
$(2.11)$, at least outwardly. In fact, we have
$$
\e{G}(z;h)=\zeta(z_1+z_2)\zeta(z_3+z_4)
\sum_{m,n\ge1}{\sigma_{z_1-z_2}(m)\sigma_{z_3-z_4}(n)\over
m^{z_1}n^{z_3}}\e{K}(m,n;\tilde{\varphi}),\;\eqno(3.5)
$$
with
$$
\tilde{\varphi}(x)={1\over\pi^2}\int_{(\alpha)}
h^{**}(s)(x/2)^{-2s}ds.
\eqno(3.6)
$$
To see this, note that $(3.1)$ implies
$$
\tilde{\varphi}(x)=\varphi(x)-\sum_{\scr{k=1}\atop\scr{k\ne6}}^7
\gamma_kJ_{2k-1}(x),\eqno(3.7)
$$
and that for $k=1,2,3,4,5,7$
$$
\eqalignno{
q_{m,n}(k)&=\sum_{l\ge1}{1\over l}S(m,n;l)J_{2k-1}\left({4\pi\over l}
\sqrt{mn}\right)\cr
&=\delta_{m,n}{(-1)^{k-1}\over2\pi},&(3.8)
}
$$
which is equivalent to the fact $\vartheta(k)=0$ for these values of
$k$ (see Lemma 2.3 of [4]). The last two identities and
$(3.3)$, $\nu=0$, give
$$
\e{K}(m,n;\tilde{\varphi})
=\e{K}(m,n;\varphi)+{\delta_{m,n}\over\pi^2}
\int_{-\infty}^\infty rh(r)\tanh\pi r\,dr\eqno(3.9)
$$
Inserting this into $(3.5)$, we recover $(2.11)$.
\par
To prove $(3.4)$, we introduce the contour $\e{L}$ which is the result
of connecting, with straight lines, the points $-\infty-ia$, $-b-ia$,
$b+ia$, $+\infty+ia$ in this order, where $a,b>0$ are arbitrary but
fixed. The analytic continuation of $h^{**}(s)$ to the region on the
left of $-i\e{L}$ is given by
$$
\eqalignno{
h^{**}(s)&={\cos\pi s\over i\pi}\int_\e{L}rh(r)\tanh\pi r\,
\Gamma(s+ir)\Gamma(s-ir)dr\cr
&-{\cos\pi s\over2i}\sum_{\scr{k=1}\atop\scr{k\ne6}}^7(-1)^k\gamma_k
\Gamma(k-\txt{1\over2}+s)\Gamma(\txt{1\over2}-k+s).&(3.10)
}
$$
Let us suppose that $s$ is on the left of $-i\e{L}$, and $t=\Im s$
tends to $+\infty$ while $-a<\sigma=\Re s$ is bounded. Since $a$
is arbitrary, this is essentially the same as to assume only that $t$
tends to
$+\infty$ while $\sigma$ remains bounded. Then an appropriate deform
of $\e{L}$ gives
$$
\eqalignno{
h^{**}(s)&={\cos(\pi s)\over i\pi}\int_{-\log t}^{\log t}
rh(r)\tanh\pi r\,\Gamma(s+ir)\Gamma(s-ir)dr\cr
&-{\cos(\pi s)\over2i}\sum_{\scr{k=1}\atop\scr{k\ne6}}^7(-1)^k\gamma_k
\Gamma(k-\txt{1\over2}+s)\Gamma(\txt{1\over2}-k+s)\cr
&+O_\sigma\left(\exp(-c(\log t)^2)\right),&(3.11)
}
$$
because of $(2.1)$. In this integral, we have, by Stirling's formula, 
$$
\eqalignno{
&\Gamma(s+ir)\Gamma(s-ir)=-2\pi it^{2\sigma-1+2it}\exp(\pi i\sigma
-2it-\pi t)\cr
&\times\left(1+t^{-1}p_1(\sigma,r^2)+t^{-2}p_2(\sigma,r^2)
+O_\sigma((|r|+1)^6 t^{-3})\right),&(3.12)
}
$$
where the polynomial $p_\nu(\sigma,\xi)$ is of degree $\nu$
in $\xi$. Thus the first term on the right of $(3.11)$ is equal to
$$
\eqalignno{
-&t^{2\sigma-1+2it}e^{-2it}\int_{-\infty}^\infty
rh(r)\tanh\pi r\cr
&\times\left(1+t^{-1}p_1(\sigma,r^2)+t^{-2}p_2(\sigma,r^2)\right)dr
+O_\sigma(t^{2\sigma-4}).&(3.13)
}
$$
By $(3.3)$, we may put this as
$$
\eqalignno{
-{\pi\over2}t^{2\sigma-1+2it}&e^{-2it}
\sum_{\scr{k=1}\atop\scr{k\ne6}}^7(-1)^k\gamma_k
\big(1+t^{-1}p_1(\sigma,(i(k-\txt{1\over2}))^2)\cr
&+t^{-2}p_2(\sigma,(i(k-\txt{1\over2}))^2)\big)
+O_\sigma(t^{2\sigma-4}),&(3.14)
}
$$
which is, by $(3.12)$, equal to
$$
{1\over2i}\sum_{\scr{k=1}\atop\scr{k\ne6}}^7(-1)^k\gamma_k
\cos \pi s\,\Gamma(k-\txt{1\over2}+s)\Gamma(\txt{1\over2}-k+s)
+O_\sigma(t^{2\sigma-4}).\eqno(3.15)
$$
We end the proof for the case $\Im s>0$. The case $\Im s<0$ is
analogous.
\bigskip
\noindent
{\bf 4.}\quad We may now return to $(3.5)$. We are about to transform 
it with Estermann's functional equation for
$$
D(s,\xi;e(a/l))=\sum_{n\ge1}n^{-s}\sigma_\xi(n)e(an/l),\quad(a,l)=1,
\eqno(4.1)
$$
basic properties of which are given in Lemma 3.7 of [4].
To this end, let us assume that
$$
\Re z_j>1-\alpha\quad (1\le j\le 4),\quad
-{5\over2}<\alpha<-{1\over2}.\eqno(4.2)
$$
Then we see readily, from $(3.5)$--$(3.6)$, that
$$
\eqalignno{
\e{G}(z;h)&={1\over\pi^2}\zeta(z_1+z_2)\zeta(z_3+z_4)\sum_{l\ge1}
{1\over l}\sum_{\scr{a\bmod l}\atop\scr{(a,l)=1}}\int_{(\alpha)}
(2\pi/l)^{-2s}h^{**}(s)\cr
&\times D(z_1+s;z_1-z_2;e(a/l))D(z_3+s,z_3-z_4;e(a^*/l))ds,&(4.3)
 }
$$
with $aa^*\equiv1\bmod l$, which converges absolutely. We consider the
following subdomain of
$(4.2)$:
$$
1-\alpha<\Re z_j<-\beta\quad (1\le j\le 4),\quad
-{7\over2}<\beta<\alpha-1.\eqno(4.4)
$$
Also we introduce the domain
$$
\Re(z_1+z_2+z_3+z_4)>3-2\beta.\eqno(4.5)
$$
We shall work temporarily in the intersection of these two domains,
which is not empty with an appropriate choice of $\alpha$, $\beta$. 
\par
We then shift the contour in $(4.3)$ to $\Re s=\beta$, which can be 
performed because of $(3.4)$ and the convexity bound for
$D(s,\xi;e(a/l))$. We encounter poles at
$s=1-z_j$ $(1\le j\le 4)$, which can be assumed, without loss of
generality, are all simple. Computing the residue, we get
$$
\e{G}(z;h)=\left(\e{R}_1+\e{G}_1\right)(z;h).\eqno(4.6)
$$
Here $\e{G}_1$ has the same expression as the right side of $(4.3)$
but with the contour $(\beta)$, and
$$
\e{R}_1(z;h)={2i\over\pi}\zeta(z_1+z_2)\zeta(z_3+z_4)
\sum_{l\ge 1}{1\over l}\sum_{\scr{a\bmod l}
\atop\scr{(a,l)=1}}\{\cdots\},\eqno(4.7)
$$
with 
$$
\eqalignno{
\{\cdots\}=&(2\pi/l)^{2(z_1-1)}\zeta(z_2-z_1+1)
l^{z_1-z_2-1}\cr
&\qquad\times D(z_3-z_1+1,z_3-z_4;e(a^*/l))h^{**}(1-z_1)\cr
+&(2\pi/l)^{2(z_2-1)}\zeta(z_1-z_2+1)
l^{z_2-z_1-1}\cr
&\qquad\times D(z_3-z_2+1,z_3-z_4;e(a^*/l))h^{**}(1-z_2)\cr
+&(2\pi/l)^{2(z_3-1)}\zeta(z_4-z_3+1)
l^{z_3-z_4-1}\cr
&\qquad\times D(z_1-z_3+1,z_1-z_2;e(a/l))h^{**}(1-z_3)\cr
+&(2\pi/l)^{2(z_4-1)}\zeta(z_3-z_4+1)
l^{z_4-z_3-1}\cr
&\qquad\times D(z_1-z_4+1,z_1-z_2;e(a/l))h^{**}(1-z_4). &(4.8)
}
$$
The $\e{R}_1(z;h)$ can be expressed in terms of the Riemann
zeta-function. To see this we assume temporarily that $z_j$ are
close to each other, while satisfying $(4.4)$--$(4.5)$ and $z_j\ne
z_{j'}$ for $j\ne j'$. Then the right side of $(4.7)$ obviously
converges to a regular function of $z$ in such a domain. With this,
let us further suppose that $\Re z_1<\Re z_3$, $\Re z_1<\Re z_4$.
We have, because of absolute convergence,
$$
\eqalignno{
&\sum_{l\ge1}l^{-z_1-z_2}\sum_{\scr{a\bmod l}
\atop\scr{(a,l)=1}}D(z_3-z_1+1,z_3-z_4;e(a/l))\cr
=&\sum_{l\ge1}l^{-z_1-z_2}\sum_{n\ge1}c_l(n)\sigma_{z_3-z_4}(n)
n^{z_1-z_3-1}\cr
=&\sum_{n\ge1}\sigma_{z_3-z_4}(n)n^{z_1-z_3-1}
\sum_{l\ge1}c_l(n)l^{-z_1-z_2}\cr
=&{1\over\zeta(z_1+z_2)}\sum_{n\ge1}\sigma_{1-z_1-z_2}(n)
\sigma_{z_3-z_4}(n)n^{z_1-z_3-1}\cr
=&{\zeta(z_3-z_1+1)\zeta(z_4-z_1+1)\zeta(z_2+z_4)\zeta(z_2+z_4)\over
\zeta(z_1+z_2)\zeta(z_3-z_1+z_3+z_4+1)},&(4.9)
}
$$
where $c_l(n)$ is the Ramanujan sum. The last expression gives a
meromorphic continuation of the initial sum to $\B{C}^4$. In this way
we find that $\e{R}_1$ exists as a meromorphic function over
$\B{C}^4$, and admits the expression
$$
\eqalignno{
\e{R}_1&(z;h)\cr
=&i{(2\pi)^{2z_1-1}h^{**}(1-z_1)
\over\pi^2\zeta(1-z_1+z_2+z_3+z_4)}\zeta(1-z_1+z_2)\zeta(1-z_1+z_3)
\zeta(1-z_1+z_4)\cr
&\qquad\times\zeta(z_2+z_3)\zeta(z_2+z_4)\zeta(z_3+z_4)\cr
+&i{(2\pi)^{2z_2-1}h^{**}(1-z_2)
\over\pi^2\zeta(1-z_2+z_1+z_3+z_4)}\zeta(1-z_2+z_1)\zeta(1-z_2+z_3)
\zeta(1-z_2+z_4)\cr
&\qquad\times\zeta(z_1+z_3)\zeta(z_1+z_4)\zeta(z_3+z_4)\cr
+&i{(2\pi)^{2z_3-1}h^{**}(1-z_3)
\over\pi^2\zeta(1-z_3+z_1+z_2+z_4)}\zeta(1-z_3+z_1)\zeta(1-z_3+z_2)
\zeta(1-z_3+z_4)\cr
&\qquad\times\zeta(z_1+z_2)\zeta(z_1+z_4)\zeta(z_2+z_4)\cr
+&i{(2\pi)^{2z_4-1}h^{**}(1-z_4)
\over\pi^2\zeta(1-z_4+z_1+z_2+z_3)}\zeta(1-z_4+z_1)\zeta(1-z_4+z_2)
\zeta(1-z_4+z_3)\cr
&\qquad\times\zeta(z_1+z_2)\zeta(z_1+z_3)\zeta(z_2+z_3).&(4.10)
}
$$ 
\par
We turn to $\e{G}_1(z;h)$; note that we assume $(4.4)$--$(4.5)$.
We replace the Estermann zeta-functions involved 
in $\e{G}_1$ by the
absolutely convergent Dirichlet series which are implied by the
functional equation. Then, after some rearrangement, we get 
$$
\e{G}_1(z;h)=\left(\e{G}_1^+
+\e{G}_1^-\right)(z;h),\eqno(4.11)
$$
where
$$
\eqalignno{
&\e{G}_1^\pm(z;h)=\zeta(z_1+z_2)\zeta(z_3+z_4)\cr
&\times\sum_{m,n\ge1}{\sigma_{z_2-z_1}(m)\sigma_{z_4-z_3}(n)\over
m^{{1\over2}(-z_1+z_2+z_3+z_4)}n^{{1\over2}(z_1+z_2-z_3+z_4)}}
\e{K}(\pm m,n;\psi_\pm).
&(4.12)
}
$$
Here
$$
\eqalignno{
\psi_+(x)&={1\over\pi^4}\int_{(\beta)}(x/2)^{2s+z_1+z_2+z_3+z_4-2}\cr
&\times\Gamma(1-z_1-s)\Gamma(1-z_2-s)\Gamma(1-z_3-s)\Gamma(1-z_4-s)\cr
&\times\Big(\cos\pi(s+\txt{1\over2}(z_1+z_2))
\cos\pi(s+\txt{1\over2}(z_3+z_4))\cr
&\qquad+\cos\txt{1\over2}\pi(z_1-z_2)
\cos\txt{1\over2}\pi(z_3-z_4)\Big)h^{**}(s)ds,
&(4.13)
}
$$
$$
\eqalignno{
\psi_-(x)&=-{1\over\pi^4}\int_{(\beta)}(x/2)^{2s+z_1+z_2+z_3+z_4-2}\cr
&\times\Gamma(1-z_1-s)\Gamma(1-z_2-s)\Gamma(1-z_3-s)\Gamma(1-z_4-s)\cr
&\times\Big(\cos\txt{1\over2}\pi(z_1-z_2)
\cos\pi(s+\txt{1\over2}(z_3+z_4))\cr
&\qquad+\cos\txt{1\over2}\pi(z_3-z_4)
\cos\pi(s+\txt{1\over2}(z_1+z_2))\Big)h^{**}(s)ds.
&(4.14)
}
$$
Invoking Weil's bound for $S(m,n;l)$, we see readily that $(3.4)$
implies that the triple sum in $(4.12)$ converges absolutely in the
domain
$$
\eqalignno{
\Big\{z:\,\Re(z_1&+z_2+z_3+z_4)>{5\over2}-2\beta,\cr
&\Re z_1,\Re z_2,\Re z_3,\Re z_4<-\beta\Big\},&(4.15)
}
$$
with any $-{7\over2}<\beta<-{5\over4}$. That is, $(4.6)$, $(4.10)$,
$(4.11)$ yield an analytic continuation of $\e{G}(z;h)$ to $(4.15)$.
Note that
$(4.15)$ contains $(4.4)$--$(4.5)$.
\bigskip
\noindent
{\bf 5.}\quad As remarked already, interesting
situations occur when we take all $z_j$ to the vicinity of $\e{X}$,
which is defined by $(2.6)$. However, the domain
$(4.15)$ does not contain such points. This necessitates analytic
continuation of
$\e{G}_1(z;h)$. It is accomplished via spectral expansions of
$\e{K}(\pm m,n;\psi_\pm)$. We shall consider first the expansion in 
the plus case, with greater details; the minus case is much simpler,
and will be treated later very briefly. The analytic continuation
itself will be established in the next section.
\smallskip
One may wish to appeal to the Kloosterman-Spectral sum formula due to
Kuznetsov, a version of which is in Theorem 2.3 of [4]. However,
our $\psi_+$ does not appear to satisfy the condition
$(2.4.6)$ postulated there, as $x$
tends to $+\infty$. Accordingly we shall instead follow the
argument employed in the proof of the theorem.
\par
To this end, let us introduce the Kloosterman-sum zeta-function:
$$
Z_{m,n}(s)=(2\pi\sqrt{mn})^{2s-1}\sum_{l\ge1}
l^{-2s}S(m,n;l).\eqno(5.1)
$$
We have, by definition,
$$
\eqalignno{
\e{K}(m,n;\psi_+)&={1\over\pi^4}\int_{(\beta)}
Z_{m,n}(s+\txt{1\over2}(z_1+z_2+z_3+z_4-1))\cr
&\times\Gamma(1-z_1-s)\Gamma(1-z_2-s)\Gamma(1-z_3-s)\Gamma(1-z_4-s)\cr
&\times\Big(\cos\pi(s+\txt{1\over2}(z_1+z_2))
\cos\pi(s+\txt{1\over2}(z_3+z_4))\cr
&\qquad+\cos\txt{1\over2}\pi(z_1-z_2)
\cos\txt{1\over2}\pi(z_3-z_4)\Big)h^{**}(s)ds.
&(5.2)
}
$$
Then we invoke Kuznetsov's spectral expansion of $Z_{m,n}(s)$ as given
in Lemma 2.5 of [4]: For any $m,n>0$ and $\Re s>{3\over4}$
$$
Z_{m,n}(s)=\left\{Z_{m,n}^{(d)}+Z_{m,n}^{(h)}+Z_{m,n}^{(c)}\right\}(s)
-{\delta_{m,n}\over2\pi}{\Gamma(s)\over\Gamma(1-s)},\eqno(5.3)
$$
where
$$
\eqalignno{
Z_{m,n}^{(d)}&={1\over2}\sin\pi s\sum_{j\ge1}\alpha_j
t_j(m)t_j(n)\Gamma(s-\txt{1\over2}
+i\kappa_j)\Gamma(s-\txt{1\over2}-i\kappa_j),\cr
Z_{m,n}^{(h)}&=\sum_{k\ge1}(2k-1)q_{m,n}(k){\Gamma(k-1+s)\over\Gamma
(k+1-s)},\cr
Z_{m,n}^{(c)}&={1\over2\pi}\sin\pi s\int_{-\infty}^\infty
{\sigma_{2ir}(m)\sigma_{2ir}(n)\over(mn)^{ir}|\zeta(1+2ir)|^2}
\Gamma(s-\txt{1\over2}+ir)
\Gamma(s-\txt{1\over2}-ir)dr.\qquad&(5.4)
}
$$
We have
$$
\e{K}(m,n;\psi_+)=\left\{\e{K}^{(d)}+
\e{K}^{(h)}+\e{K}^{(c)}+\e{K}^{(\delta)}\right\}(m,n;\psi_+),
\eqno(5.5)
$$
with an obvious arrangement of terms.
\par
We claim that
$$
\e{K}^{(d)}(m,n;\psi_+)=\sum_{j\ge1}\alpha_j
t_j(m)t_j(n)\Psi_+(z;\kappa_j),\eqno(5.6)
$$
where
$$
\eqalignno{
&\Psi_+(z;r)={1\over2\pi^4}\int_{-i\infty}^{i\infty}
h^{**}(s)\Gamma(s-1+\txt{1\over2}(z_1+z_2+z_3+z_4)+ir)\cr
&\times\Gamma(s-1+\txt{1\over2}(z_1+z_2+z_3+z_4)-ir)
\Gamma(1-z_1-s)\Gamma(1-z_2-s)
\cr&\times\Gamma(1-z_3-s)\Gamma(1-z_4-s)
\sin\pi(s+\txt{1\over2}(z_1+z_2+z_3+z_4-1))\cr
&\times\Big(\cos\pi(s+\txt{1\over2}(z_1+z_2))
\cos\pi(s+\txt{1\over2}(z_3+z_4))\cr
&\hskip2cm+\cos\txt{1\over2}\pi(z_1-z_2)
\cos\txt{1\over2}\pi(z_3-z_4)\Big)ds.&(5.7)
}
$$
Here the path of integral is curved to ensure that the poles of the
first three factors of the integrand lie to the left, and those of
other factors to the right, respectively, while parameters are such
that the path can be drawn. 
\par
This is easy to confirm. In fact, with $r\in\B{R}$ and
$z$ in $(4.15)$, we may take the vertical line $(\beta)$ as the path.
Then $(3.4)$ and the Stirling formula implies readily
that the integral of the absolute value of the integrand in $(5.7)$
is $\ll_z(|r|+1)^{-4}$. Also we have the bounds $(1.4)$ and
$$
\sum_{\kappa_j\le K}\alpha_j\ll K^2\eqno(5.8)
$$
(see $(2.3.2)$ of [4]). By these facts,  we see that $(5.6)$ converges
uniformly and absolutely in any compactum of $(4.15)$. Obviously we
have
$$
\e{K}^{(c)}(m,n;\psi_+)={1\over\pi}\int_{-\infty}^\infty
{\sigma_{2ir}(m)\sigma_{2ir}(n)\over(mn)^{ir}|\zeta(1+2ir)|^2}
\Psi_+(z;r)dr,\eqno(5.9)
$$
as well. Incidentally we have proved that for
$z$ in $(4.15)$ and $r\in\B{R}$
$$
\Psi_+(z;r)\ll (|r|+1)^{-4}.\eqno(5.10)
$$
\par
Further, we have
$$
\e{K}^{(h)}(m,n;\psi_+)={2\over\pi}\sum_{k\ge1}(-1)^k
(2k-1)q_{m,n}(k)\Psi_+(z;i(\txt{1\over2}-k)).\eqno(5.11)
$$
This is a simple consequence of
$$
q_{m,n}(k)\ll{(2\pi\sqrt{mn})^{2k-1}\over\Gamma(2k)}
\sum_{l\ge1}|S(m,n;l)|l^{-2}\eqno(5.12)
$$
(see $(2.2.7)$ of [4]) and, for $z$ in $(4.15)$,
$$
\eqalignno{
\Psi_+(z;&i(\txt{1\over2}-k))
={(-1)^k\over2\pi^3}\int_{(\beta)}h^{**}(s)
{\Gamma(s+{1\over2}(z_1+z_2+z_3+z_4-3)+k)
\over\Gamma(k-s-{1\over2}(z_1+z_2+
z_3+z_4-3))}\cr
&\times\Gamma(1-z_1-s)\Gamma(1-z_2-s)
\Gamma(1-z_3-s)\Gamma(1-z_4-s)\cr
&\times\Big(\cos\pi(s+\txt{1\over2}(z_1+z_2))
\cos\pi(s+\txt{1\over2}(z_3+z_4))\cr
&\hskip2cm+\cos\txt{1\over2}\pi(z_1-z_2)
\cos\txt{1\over2}\pi(z_3-z_4)\Big)ds.&(5.13)
}
$$
\par
Corresponding to $(5.10)$ we have
$$
\Psi_+(z;i(\txt{1\over2}-k))\ll k^{-4}\quad(k=1,2,3\ldots)\eqno(5.14)
$$
uniformly for any compactum of
$$
\eqalignno{
\Big\{z:&\, {5\over2}-2\beta<\Re(z_1+z_2+z_3+z_4)<6,\cr
&\Re z_1,\Re z_2,\Re z_3,\Re z_4<-\beta\Big\},
\quad -{7\over4}<\beta<-{5\over4}.&(5.15)
}
$$
To show this, we move the contour in $(5.13)$ to the one which is
the result of connecting, with straight lines, the points
$\beta-i\infty$, $\beta-{1\over2}ki$, $-4-\txt{1\over2}ki$,
$-4+\txt{1\over2}ki$,
$\beta+\txt{1\over2}ki$, and $\beta+i\infty$, in this order. We
encounter only one singularity which is the simple pole of
$h^{**}(s)$ at $s=-{7\over2}$. The residue is $\ll
k^{-10+\Re(z_1+z_2+z_3+z_4)}$. The integral  over the new path is
readily estimated to be $O(k^{-4})$ with $(3.4)$ and Stirling's
formula. Thus we have $(5.14)$.
\par
With this, we rearrange $(5.11)$, following the argument on p.\ 66 
of [4], and obtain, for $z$ in $(5.15)$,
$$
\left(\e{K}^{(h)}+\e{K}^{(\delta)}\right)(m,n;\psi_+)
=\sum_{k\ge1}\sum_{j=1}^{\vartheta(k)}
\alpha_{j,k}t_{j,k}(m)t_{j,k}(n)\Psi_+(z;i(\txt{1\over2}-k)),
\eqno(5.16)
$$
with $\alpha_{j,k}=16\Gamma(2k)(4\pi)^{-2k-1}
|\varrho_{j,k}(1)|^2$. Note that we have
$$
\sum_{j=1}^{\vartheta(k)}\alpha_{j,k}\ll k\eqno(5.17)
$$
(see $(2.2.10)$ of [4]).
\smallskip
Collecting $(4.12)^+$, $(5.5)$, $(5.6)$, $(5.9)$ and $(5.16)$, we
find that for $z$ in $(5.15)$
$$
\e{G}_1^+(z;h)=
\left\{\e{H}_++\e{H}_0+\e{C}\right\}(z^*;\Psi_+(z;\cdot)),
\eqno(5.18)
$$
where for $z=(z_1,z_2,z_3,z_3)$ 
$$
z^*=(z_1^*,z_2^*,z_3^*,z_4^*),\quad
z_j^*={1\over2}(z_1+z_2+z_3+z_4)-z_j;\eqno(5.19)
$$ 
$\e{H}_+$, $\e{C}$ are as in $(2.3)$, but with an obvious extension
of notation, and
$$
\e{H}_0(z;h)
=\sum_{k\ge1}\sum_{j=1}^{\vartheta(k)}
\alpha_{j,k}H_{j,k}(z_1)H_{j,k}(z_2)H_{j,k}(z_3)H_{j,k}(z_4)
h(i(\txt{1\over2}-k)).\eqno(5.20)
$$
The absolute convergence that is needed to carry out this rearranging
procedure is guaranteed by the uniform bounds $(1.4)$, $(1.8)$ for
Hecke eigenvalues, and by the bounds $(5.8)$, $(5.10)$,
$(5.14)$, $(5.17)$.
\smallskip
We turn to $\e{G}_1^-(z;h)$, or
the spectral decomposition of $\e{K}(-m,n;\psi_-)$. Note that $z$
is now to be in the domain $(5.15)$. This time we appeal to Theorem
2.5 of [4]. The condition $(2.4.6)$ of [4] is easily seen to be
satisfied by
$\psi_-(x)$, because of the exponential decay of the integrand in
$(4.14)$: According as $x\to+0$ and $+\infty$, we use the contours
$(\beta)$ and $(-3)$, respectively. 
\par
Skipping the details, which are
easy to fill, we have, corresponding to $\Psi_+$,
$$
\eqalignno{
&\Psi_-(z;r)=-{\cosh\pi r\over2\pi^4}\int_{-i\infty}^{i\infty}
h^{**}(s)\Gamma(s-1+\txt{1\over2}(z_1+z_2+z_3+z_4)+ir)\cr
&\times\Gamma(s-1+\txt{1\over2}(z_1+z_2+z_3+z_4)-ir)
\Gamma(1-z_1-s)\Gamma(1-z_2-s)
\cr&\times\Gamma(1-z_3-s)\Gamma(1-z_4-s)
\Big(\cos\txt{1\over2}\pi(z_1-z_2)
\cos\pi(s+\txt{1\over2}(z_3+z_4))\cr
&\hskip2cm+\cos\pi(s+\txt{1\over2}(z_1+z_2))
\cos\txt{1\over2}\pi(z_3-z_4)\Big)ds,&(5.21)
}
$$
with the same path as in $(5.7)$. Also we have, uniformly for
bounded $z$ in $(5.15)$,
$$
\Psi_-(z;r)\ll(|r|+1)^{-4}.\eqno(5.22)
$$
With this, we have, in much the same way as before,
$$
\e{G}_1^-(z;h)=
\left\{\e{H}_-+\e{C}\right\}(z^\ast;\Psi_-(z;\cdot)).\eqno(5.23)
$$
\smallskip
Collecting $(4.6)$, $(4.11)$, $(5.18)$, $(5.23)$, we obtain, in
the domain $(5.15)$,
$$
\eqalignno{
\e{G}(z;h)&=\e{R}_1(z;h)+\e{H}_+(z^*;\Psi_+(z;\cdot))
+\e{H}_-(z^*;\Psi_-(z;\cdot))\cr
&+\e{H}_0(z^*;\Psi_+(z;\cdot))
+\e{C}(z^*;(\Psi_++\Psi_-)(z;\cdot)).&(5.24)
}
$$
\bigskip
\noindent
{\bf 6.}\quad We shall continue analytically the last expression for
$\e{G}(z;h)$. The article [2] lacks this decisive procedure
altogether.
\smallskip
To this end, we introduce
$$
\eqalignno{
\Xi_1&(z;r)=\int_{-i\infty}^{i\infty}h^{**}(s)
{\Gamma\left(ir+{1\over2}(z_1+z_2+z_3+z_4)-1+s\right)\over
\Gamma\left(ir-{1\over2}(z_1+z_2+z_3+z_4)+2-s\right)}\cr
&\times\Gamma(1-z_1-s)\Gamma(1-z_2-s)\Gamma(1-z_3-s)
\Gamma(1-z_4-s)\cr
&\hskip1cm\times\cos\pi\left(s+\txt{1\over2}(z_1+z_2)\right)
\cos\pi\left(s+\txt{1\over2}(z_3+z_4)\right)ds,&(6.1)\cr
\Xi_2&(z;r)=\int_{-i\infty}^{i\infty}h^{**}(s)
{\Gamma\left(ir+{1\over2}(z_1+z_2+z_3+z_4)-1+s\right)\over
\Gamma\left(ir-{1\over2}(z_1+z_2+z_3+z_4)+2-s\right)}\cr
&\times\Gamma(1-z_1-s)\Gamma(1-z_2-s)\Gamma(1-z_3-s)
\Gamma(1-z_4-s)ds,&(6.2)\cr
}
$$
where the paths are the same and separate the poles of 
$h^{**}(s)\Gamma(ir+{1\over2}(z_1+z_2+z_3+z_4)-1+s)$ to the left and
those of $\Gamma(1-z_1-s)\Gamma(1-z_2-s)\Gamma(1-z_3-s)
\Gamma(1-z_4-s)$ to the right; the $(z;r)\in\B{C}^5$ is assumed to be
such that the paths can be drawn. The integrand of $(6.2)$ is of
exponential decay, and the bound
$(3.4)$ assures amply the convergence of $(6.1)$. Thus $\Xi_1(z;r)$
and $\Xi_2(z;r)$ are regular for those $(z;r)$ indicated.
Shifting the paths appropriately we see also that they are
meromorphic over $\B{C}^5$; the polar divisors arise from those
$(z;r)$ with which the path cannot be drawn. 
\par
In particular, if $z$
is in $(5.15)$ and $r\in\B{R}\cup
\{i({1\over2}-k):k=1,2,3,\ldots\}$, then the line $(\beta)$ can
be used in $(6.1)$ and $(6.2)$. With this observation, 
we find, after some computation, that
$$
\eqalignno{
\Psi_+(z;r)&=i{\Xi_1(z;r)-\Xi_1(z;-r)\over4\pi^3\sinh\pi
r}\cr &+i{\cos{1\over2}\pi(z_1-z_2)
\cos{1\over2}\pi(z_3-z_4)
\over4\pi^3\sinh\pi
r}\left\{\Xi_2(z;r)-\Xi_2(z;-r)\right\},\qquad&(6.3)\cr
\Psi_-(z;r)&={i\over4\pi^3\sinh\pi
r}\Big(\cos\txt{1\over2}\pi(z_1-z_2)
\cos\pi\left(\txt{1\over2}(z_1+z_2)+ir\right)\cr
&+\cos\txt{1\over2}\pi(z_3-z_4)
\cos\pi\left(\txt{1\over2}(z_3+z_4)+ir\right)\Big)
\Xi_2(z;r)\cr
&-{i\over4\pi^3\sinh\pi
r}\Big(\cos\txt{1\over2}\pi(z_1-z_2)
\cos\pi\left(\txt{1\over2}(z_1+z_2)-ir\right)\cr
&+\cos\txt{1\over2}\pi(z_3-z_4)
\cos\pi\left(\txt{1\over2}(z_3+z_4)-ir\right)\Big)
\Xi_2(z;-r).&(6.4)
}
$$
These identities obviously yield meromorphic continuations of
$\Psi_\pm(z;r)$ to $\B{C}^5$; that is, the two functions exist as far
as the respective right sides are finite.   Note that we have, for
integers $k\ge1$,
$$
\eqalignno{
\Psi_+(&z;i(\txt{1\over2}-k))={(-1)^k\over2\pi^3}
\Xi_1(z;i(\txt{1\over2}-k))\cr
&+{(-1)^k\over2\pi^3}\cos\txt{1\over2}\pi(z_1-z_2)
\cos\txt{1\over2}\pi(z_3-z_4)
\Xi_2(z;i(\txt{1\over2}-k)),&(6.5)
}
$$
provided the right side is finite.  Compare the above with Lemma
4.4 of [4].
\par
We shall next estimate $\Xi_1$ and $\Xi_2$ as functions of
$r$. To this end we introduce the set
$$
\left\{\hbox{$z\in\B{C}^4:\,$ none of $z_j$ are equal to 
$\displaystyle{9\over2}+a$
with an integer $a\ge0 $}\right\}.\eqno(6.6)
$$
Note that this contains $(5.15)$. 
Both $\Xi_1$ and $\Xi_2$ are
regular in any compactum of
$(6.6)$, provided $|\Re r|$ is sufficiently large. Moreover, if $r$
tends to infinity in any fixed horizontal strip, we have
$$
\eqalignno{
\Xi_1(z;r)&\ll |r|^{-4}e^{\pi|r|},&(6.7)\cr
\Xi_2(z;r)&\ll |r|^{-10+\Re(z_1+z_2+z_3+z_4)},&(6.8)
}
$$
with implied constants depending only on 
$\max(|z_1|, |z_2|, |z_3|, |z_4|, |\Im r|)$. In fact, in this
situation the paths in $(6.1)$ and $(6.2)$ can be drawn, and the
regularity assertion is immediate. To prove the last bounds as $\Re r$
is positive and large, we move the path to the one which is the
result of connecting, with straight lines, the points $M-i\infty$,
$M-{1\over2}i|r|$, $-M-{1\over2}i|r|$,
$-M+i\infty$ in this order, where $M$ is an
integer such that $|r|>2M>8\max(|z_1|, |z_2|, |z_3|, |z_4|, |\Im
r|)$. The singularities we encounter are only poles of $h^{**}(s)$.
The residues contribute $\ll_M|r|^{-10+\Re(z_1+z_2+z_3+z_4)}$ to both
$\Xi_1$ and $\Xi_2$. To estimate the integrals on the new path, we
note that by $(3.4)$ and Stirling's formula the integrand of $\Xi_1$
is estimated to be
$$
\eqalignno{
&\ll_M|s+ir|^{\Re(s+{1\over2}(z_1+z_2+z_3+z_4)+ir)-{3\over2}}\cr
&\times|s-ir|^{\Re(s+{1\over2}(z_1+z_2+z_3+z_4)-ir)-{3\over2}}
|s|^{-\Re(2s+z_1+z_2+z_3+z_4)-2}e^{\pi|r|},&(6.9)
}
$$
and that of $\Xi_2$ by the same expression but with $e^{-\pi|s|}$
in place of $e^{\pi|r|}$. The bounds $(6.7)$--$(6.8)$ follow
immediately.
\par
Hence we have, uniformly for any compactum in $(6.6)$,
$$
\eqalignno{
\Psi_+(z;r)&\ll |r|^{-4},&(6.10)\cr
\Psi_-(z;r)&\ll |r|^{-10+\Re(z_1+z_2+z_3+z_4)},&(6.11)
}
$$
as $r$ tends to infinity in any fixed horizontal strip. In much the
same way one may show that uniformly for all integers $k\ge1$ and
for any compactum in $(6.6)$
$$
\Psi_+(z;i(\txt{1\over2}-k))\ll k^{-4}+k^{-10+\Re(z_1+z_2+z_3+z_4)}.
\eqno(6.12)
$$
\par
Now, the combination $(6.10)$--$(6.12)$ implies an important
assertion: the spectral sums 
$\e{H}_\pm(z^*;\Psi_\pm(z;\cdot))$ and 
$\e{H}_0(z^*;\Psi_+(z;\cdot))$ are all meromorphic in the domain
$$
\prod_{j=1}^4\left\{z_j:\,{1\over3}<\Re z_j<{3\over2}\right\}.
\eqno(6.13)
$$
To show this, we invoke the following large sieve estimates: For any
fixed $\varepsilon>0$
$$
\eqalignno{
\sum_{\kappa_j\le K}\alpha_j&|H_j(s)|^4,\;\sum_{k\le K}
\sum_{j=1}^{\vartheta(k)}\alpha_{j,k}|H_{j,k}(s)|^4\cr
&\ll_{\varepsilon,s} K^{2+8\max(0,{1\over2}-\Re s)+\varepsilon}.
&(6.14)
}
$$
Namely, it is enough to prove that for $z$ in $(6.13)$
$$
\sum_{j=1}^4\max\left(0,\txt{1\over2}-\Re z_j^*\right)<1,\eqno(6.15)
$$
which is, however, elementary.
\par
We consider the contribution of the continuous spectrum in $(5.24)$.
Let $z$ be in $(5.15)$. By $(6.3)$--$(6.4)$ we have
$$
\e{C}(z^*;(\Psi_++\Psi_-)(z;\cdot))=\e{Y}(z;h),
\eqno(6.16)
$$ 
where
$$
\eqalignno{
&\e{Y}(z;h)={1\over\pi}\int_{(0)}
\zeta(z_1^*+\xi)\zeta(z_1^*-\xi)\zeta(z_2^*+\xi)\zeta(z_2^*-\xi)\cr
&\times\zeta(z_3^*+\xi)\zeta(z_3^*-\xi)\zeta(z_4^*+\xi)
\zeta(z_4^*-\xi)
{\Theta(z;-i\xi)\over\zeta(1+2\xi)\zeta(1-2\xi)}d\xi,&(6.17)
}
$$
with
$$
\eqalignno{
\Theta(z;r)={1\over2\pi^3\sinh\pi r}&\Big[\Xi_1(z;r)+
\Big\{\cos\pi(ir+\txt{1\over2}(z_1+z_2))
\cos\txt{1\over2}\pi(z_1-z_2)\cr
&+\cos\pi(ir+\txt{1\over2}(z_3+z_4))
\cos\txt{1\over2}\pi(z_3-z_4)\cr
&+\cos\txt{1\over2}\pi(z_1-z_2)
\cos\txt{1\over2}\pi(z_3-z_4)\Big\}
\Xi_2(z;r)\Big].&(6.18)
}
$$
Observe first that if $z$ is in
$(5.15)$ then $\Xi_1(z;-i\xi)$ and
$\Xi_2(z;-i\xi)$ are regular for $\Re\xi\ge -{5\over4}$. In fact, on
noting the remark made after $(6.6)$, possible singularities of these
functions of $\xi$ are at
$$
-z_1^*-a_1,\;-z_2^*-a_2,\;-z_3^*-a_3,\;-z_4^*-a_4,\eqno(6.19)
$$
with non-negative integers $a_j$, and the real parts of these points
are all less than $-{5\over4}$. We move the contour of the
integral for $\e{Y}$ to the one that is the result of connecting,
with straight lines, the points $-i\infty$, $-Ni$,
$[N]+{3\over4}-Ni$, $[N]+{3\over4}+Ni$, $Ni$, $+i\infty$, in this
order. Here $N>0$ is large and such that $\zeta(s)\ne0$ for $\Im
s=2N$. Assuming further that $|z_j|<{1\over4}N$ $(1\le j\le4)$, we
have
$$
\eqalignno{
\e{Y}(z;h)=&[\hbox{the new integral}]
\cr
+2i\times &[\hbox{the sum of
residues at $z_j^*-1$  $(1\le j\le4)$}]\cr
+2i\times &[\hbox{the sum of
residues at ${1\over2}\rho$ with $|\Im\rho|<2N$}],&(6.20)
}
$$
where $\rho$ is a complex zero of $\zeta(s)$. Note that we have
applied the functional equation for the zeta-function to
$\zeta(1-2\xi)\sin\pi\xi$, and also that we have assumed that those
relevant singularities are all simple poles, which obviously does not
cause any loss of generality. This expression yields a meromorphic
continuation of
$\e{Y}$ to the domain $(6.13)$. In fact, the last two terms in
$(6.20)$ are meromorphic over $\B{C}^4$; and the uniform convergence
of the new integral for any compactum in $(6.13)$ is a consequence of
$(6.7)$--$(6.8)$ and the easy bound
$$
\int_0^T|\zeta(\sigma+it)|^8dt\ll T^{2+8\max(0,{1\over2}-\sigma)}. 
\eqno(6.21)
$$
Thus $\e{C}(z^*;(\Psi_++\Psi_-)(z;\cdot))$ exists as a meromorphic
function in $(6.13)$. 
\smallskip
Observe that $(5.15)$ and $(6.13)$ are not disjoint. Namely, we have
established the desired analytic continuation of the identity
$(5.24)$.
\bigskip
\noindent
{\bf 7.}\quad We now specialize the above by taking $z$ close to
the set $\e{X}$. The functions $\Xi_1(z;r)$ and
$\Xi_2(z;r)$ are regular in such $z$, for each $r\in\B{R}$; one may
use, for instance, the line $({1\over4})$ as the path in $(6.1)$ and
$(6.2)$. This implies, via $(6.3)$--$(6.4)$, $\Psi_\pm(z;r)$ are
regular in the same way. Similarly, $\Psi_+(z;i({1\over2}-k))$ is
regular. Hence
$\e{H}_\pm(z^*;\Psi_\pm(z;\cdot))$ and $\e{H}_0(z^*;\Psi_+(z;\cdot))$
are all regular for the current specialization of $z$.
As to $\e{Y}(z;h)$, we shift the contour on the right of $(6.20)$
back to the imaginary axis. We have
immediately
$$
\e{C}(z^*;(\Psi_++\Psi_-)(z;\cdot))
=\left(\e{C}_0+\e{R}_0\right)(z^*;
(\Psi_++\Psi_-)(z;\cdot)).\eqno(7.1)
$$
\smallskip
Collecting $(2.7)$, $(5.24)$, and $(7.1)$, we obtain:
\smallskip
\noindent
{\bf Theorem.}\quad{\it Let $z$ be either in $\e{X}$ or close to it.
Then we have the functional equation
$$
\eqalignno{
&\left(\e{H}_++\e{C}_0+\e{R}_0-\e{R}_1\right)(z;h)\cr
&=\e{H}_+(z^*;\Psi_+(z;\cdot))+\e{H}_-(z^*;\Psi_-(z;\cdot))\cr
&+\e{H}_0(z^*;\Psi_+(z;\cdot))
+(\e{C}_0+\e{R}_0)(z^*;(\Psi_++\Psi_-)(z;\cdot)),&(7.2)
}
$$
with the conventions introduced above.
}
\smallskip
\noindent
This is a corrected version of Theorem 15 of [2].
\bigskip
\noindent
{\csc Concluding Remark.}\quad The device shown in Section 3 was
sketched in Kuz\-netsov's letter (received on May
16, 1991), and was attributed to him in [3]  by the present author,
to whom the above proof is due.  It
should, however, be recorded, to avoid any confusion in the
future, that a few variants of the device had been known to
specialists since the publication of Murty's idea [5]. 
\par 
As a matter of fact, the device is by no means
of absolute necessity in treating the spectral fourth moment of Hecke
$L$-functions. For instance, in [1] it is replaced, though implicitly,
by a series of truncation procedures. What is really indispensable is
the perfection of analytic continuation, which in the present context
is achieved only in Section 6. This aspect
becomes apparent if one considers $\e{H}_-(z;h)$ instead. There the
device is irrelevant, and the whole argument is clearly built upon
the analytic continuation of the spectral decomposition. It should
perhaps be added that in [1] the continuation procedure is hidden in
a use of the explicit formula for the binary additive divisor sum.

\vskip1cm
\centerline{\bf REFERENCES}
\medskip
\item{[1]} M. Jutila and Y. Motohashi. A note on the mean value of the
zeta and $L$-functions.\ XI. Proc.\ Japan Acad., {\bf78A}, 1--6
(2002).
\item{[2]} N.V. Kuznetsov. Sums of Kloosterman sums and the eighth
power moment of the Riemann zeta-function. T.I.F.R. Stud.\ Math.,
{\bf12}, 57--117 (1989).
\item{[3]} Y. Motohashi. Kuznetsov's paper on the eighth power moment
of the Riemann zeta-function (Revised).\ Part I.
Manuscript, June 22, 1991.
\item{[4]} ---. {\it Spectral Theory of the
Riemann Zeta-Function\/}. Cambridge Univ.\ Press,
Cambridge, 1997.
\item{[5]} R. Murty. Oscillations of Fourier coefficients of modular
forms. Math.\ Ann., {\bf 262}, 431--446 (1983).
\vskip 1cm
\font\small=cmr8
{\small\noindent 
Yoichi Motohashi
\par\noindent
Honkomagome 5-67-1-901, Tokyo 113-0021, Japan
\par\noindent
Email: am8y-mths@asahi-net.or.jp 
}

\bye